\documentclass[11pt]{amsart}
\usepackage{amssymb, amsmath,mathrsfs}
\usepackage{epsfig}

\newtheorem{thm}{Theorem}
\newtheorem{cor}[thm]{Corollary}
\newtheorem{prop}[thm]{Proposition}
\newtheorem{lem}[thm]{Lemma}
\newtheorem{exmp}[thm]{Example}
\theoremstyle{definition}
\newtheorem{defn}[thm]{Definition}

\newcommand{\ol}{\overline}
\newcommand{\sq}{\subseteq}

\newcommand{\sm}{\smallsetminus}

\newcommand{\les}{\leqslant}

\newcommand{\gd}{$G_{\delta}$-}

\newcommand{\Om}{\Omega}

\newcommand{\mc}{\mathcal{C}}

\newcommand{\n}{\varnothing}

 \def\T{{\mathcal T}} 
  
 \let\alp\alpha \let\bt\beta  
\let\dlt\delta   \let\kp\kappa  
  \let\w\omega \let\phi\varphi 
\let\<\langle \let\>\rangle \let\nowt\varnothing
 \def\ol{\overline}

\newcommand{\mn}{monotonically normal\:}
\newcommand{\mdn}{\text{m}\delta\text{n}}
\newcommand{\mddn}{\text{m}\delta\dlt\text{n}}

\newcommand{\cdn}{c\(\delta\)n}
\newcommand{\cdns}{c\(\delta\)n\:}
\newcommand{\mcdn}{mc\(\delta\)n}
\newcommand{\mcdns}{mc\(\delta\)n\:}
\newcommand{\wcdn}{wc\(\delta\)n}
\newcommand{\wcdns}{wc\(\delta\)n\:}

\newcommand{\lw}{\mathbb L_{\w_1}}
\newcommand{\lww}{\mathbb L_{\w_2}}

\begin{document}

\setlength{\parindent}{0in}
\setlength{\parskip}{2ex}

\title{Monotone versions of $\delta$-normality}
\author{Chris Good} 
\email{C.Good@bham.ac.uk}
\author{Lylah Haynes}
\email{HaynesL@maths.bham.ac.uk}
\address{School of Mathematics, University of Birmingham, Birmingham B15
2TT, UK}

\begin{abstract}
We continue the study of properties related to monotone countable
paracompactness, investigating various monotone versions of
\(\delta\)-normality. We factorize monotone normality and
stratifiability in terms of these weaker properties.
\end{abstract}
\maketitle



\section{Introduction} Dowker \cite{d} proves that the product of
a space $X$ and the closed unit interval $[0,1]$ is normal iff $X$
is both normal and countably paracompact. Mack \cite{mack2} proves
that a space \(X\) is countably paracompact iff \(X \times [0,1]\)
is \(\delta\)-normal and that every countably paracompact space is
$\dlt$-normal (see below for definitions).

In \cite{gks} and its sequels \cite{gy,gk}, the first author
\textit{et al.} introduce and study a monotone version of countable
paracompactness (MCP) closely related to stratifiabilty.  In
\cite{gh}, the current authors consider various other possible
monotone versions of countable paracompactness and the notion of
$\mdn$ (monotone $\delta$-normality) arises naturally in this study.
It turns out that MCP and $\mdn$ are distinct properties and that,
if \(X \times [0,1]\) is m\(\delta\)n, then \(X\) (and hence
$X\times[0,1]$) is MCP.

In this paper we take a closer look at monotone
versions of $\dlt$-normality.

Our notation and terminology are standard as found in \cite{eng} or
\cite{kunen}. All spaces are assumed to be \(T_1\) and regular.

\section{Monotone versions of $\dlt$-normality}

\begin{defn} Let $X$ be a space.
A subset \(D\) of $X$ is said to be a regular \gd set iff there
exist open sets \(U_n\), $n\in\w$, 
such that \(D \sq U_n\) for each \(n\) and 
\(D=\bigcap_{n \in \omega} \ol{U}_n\).
\end{defn}

Clearly, a set \(D\) is a regular \gd set iff there exist open sets \(U_n\), \(n \in \omega\), such that \(D = \bigcap_{n \in \omega} {U}_n =\bigcap_{n \in \omega} \ol{U}_n\).

\begin{defn}
$X$ is said to be \emph{$\dlt$-normal} \cite{mack2} iff any two
disjoint closed sets, one of which is a regular \gd set, can be
separated by open sets.

$X$ is said to be \emph{weakly $\dlt$-normal} \cite{sigma} iff any
two disjoint regular \gd sets can be separated by open sets.
\end{defn}

We note in passing the following facts about regular \gd sets.
Finite unions and countable intersections of regular \gd sets
are again regular $G_\dlt$. If \(X\) is \(T_3\), for every \(x \in X\) and every
open neighbourhood \(V\) of
\(x\) there exists a regular \gd set \(K\) such that \(x \in K \sq V\).
In any space \(X\), the zero-sets are regular \gd sets and
so in a normal space \(X\),
if \(C\) is a closed set contained in an open set
\(U\), then there exists an open set \(W\) such that
\(W\) is the complement of
a regular \gd set and \(C \sq W \sq \ol{W} \sq U\).
If \(E\) is a regular \gd set in \(X\), then \(E \times \{\alpha\}\) is a
regular \gd set in \(X \times M\) for any infinite compact metrizable space
\(M\) and \(\alpha \in M\).
If \(Y\) is any compact space, since the projection map is both closed and open,
then the projection of a regular \gd set in \(X
\times Y\) is itself a regular \gd set in \(X\).
On the other hand, a regular \gd subset of a regular
\gd subset of $X$ is not necessarily a regular \gd set in $X$:
for example, the \(x\)-axis, $A$, is a regular \gd subset of the
Moore plane and every subset of $A$ is a regular \gd subset in $A$.

Let us make the following definition.
\begin{defn}
Let $X$ be a space and $\mc$ be a collection of pairs of disjoint
closed sets. We shall say that $H$ is a $\mc$-mn operator on $X$ iff
$H$ assigns to each pair $(C,D)\in\mc$ an open set $H(C,D)$ such that
\begin{enumerate}
\item $C\subseteq H(C,D)\subseteq \ol{H(C,D)}\subseteq X\setminus D$,
\item if $C\subseteq C'$ and $D'\subseteq D$, then
$H(C,D)\subseteq H(C',D')$.
\end{enumerate}
\end{defn}

\begin{defn} Let $H$ be a $\mc$-mn operator on $X$.
\begin{enumerate}
\item If $\mc$ is the collection of pairs of
disjoint closed subsets of $X$, then $X$ is \mn.
\item \label{deflmdn} If $\mc$ is the collection of disjoint closed subsets
$(C,D)$ such that $C$ is a regular $G_\dlt$-set,
then $X$ is left monotonically $\dlt$-normal or l$\mdn$.
\item If
$\mc$ is the collection of pairs of
disjoint closed subsets of $X$ at least one of which
is a regular $G_\dlt$-set, then $X$ is monotonically
$\dlt$-normal or $\mdn$.
\item If
$\mc$ is the collection of pairs of
disjoint regular \gd subsets of $X$, then $X$ is $\mddn$.
\end{enumerate}
\end{defn}

It can easily be shown that right monotone $\dlt$-normality
(where $D$, rather than $C$,  is assumed to be a regular \gd set)
is equivalent to lm\(\delta\)n.

Note that, replacing $H(C,D)$ with $H(C,D)\sm \ol{H(D,C)}$ if
necessary, we may assume that $H(C,D)\cap H(D,C)=\nowt$ whenever $H$
is an mn, $\mdn$ or $\mddn$ operator.

There are a number of characterizations of monotone normality,
amongst them the equivalence of conditions (1) and (2) in Theorem
\ref{mngr2} (see \cite{g}) (the proof of the extension stated here
is routine). Mimicking the proof of this characterization, we obtain
the hierarchy of monotone versions of \(\delta\)-normality listed in
Theorem \ref{ul1}.

\begin{thm} \label{mngr2}
The following are equivalent for a space \(X\):
\newcounter{romcount}
\begin{enumerate}
\item \(X\) is monotonically normal.
\item \label{part2} There is an operator \(\psi\) assigning to each open set \(U\) in
\(X\) and \(x \in U\), an open set \(\psi(x,U)\) such that
\begin{enumerate}
\item \(x \in \psi(x,U)\),
\item if \(\psi(x,U) \cap \psi(y,V) \neq \n\), then either \(x \in
V\) or \(y \in U\).
\end{enumerate}
\item  There is an operator \(\psi\) as in (2) such that, in addition,
\(\psi(x,U) \sq U\).
\item  There is an operator \(\psi\) as in (2) such that, in addition,
\(\ol{\psi(x,U)} \sq U\).
\end{enumerate}
\end{thm}

In Theorem \ref{mngr2}, monotone normality is characterized in terms
of an operator assigning an open set to each point \(x\) and open neighbourhood $U$
of \(x\).
We define several new properties, analogous to these
characterizations, by considering an operator acting on a
regular \gd set \(L\) and an open set containing \(L\).

\begin{defn} \label{defnuln}
A space \(X\) is \textit{weakly coherently} \(\delta\)\textit{-normal (\wcdn)} iff there is an operator
\(\phi\) assigning to each regular \gd set \(L\) and open set \(U\) containing
\(L\), an open set \(\phi(L,U)\) such that
\begin{enumerate}
\item \(L \sq \phi(L,U)\),
\item if \(\phi(L,U) \cap \phi(K,V) \neq \n\) then either \(L \cap V
\neq \n\) or \(K \cap U \neq \n\).
\end{enumerate}
\(X\) is \textit{coherently \(\delta\)-normal (\cdn)} if in addition,
\begin{enumerate}
\item[(3)] \(L\subseteq \phi(L,U)\sq \ol{\phi(L,U)} \sq U\).
\end{enumerate}
\(X\) is \textit{monotonically coherently \(\delta\)-normal (\mcdn)} if in addition,
\begin{enumerate}
\item[(4)] if \(L \sq L'\) and \(U \sq U'\) then \(\phi(L,U) \sq \phi(L',U')\).
\end{enumerate}
\end{defn}

If $\phi$ is an operator witnessing that $X$ is \wcdn, there
is no assumption that $\phi(L,U)$ is monotone in $L$ or $U$ nor
that it is a subset of $U$. We have the following proposition.

\begin{prop} \label{monul}
Suppose that \(X\) is  \wcdn. Then there is a \wcdns
operator $\phi$ on $X$ such that:
\begin{enumerate}
\item \(L \sq \phi(L,U) \sq U\) and
\item if \(L \sq L'\) and \(U \sq U'\), then \(\phi(L,U) \sq\phi(L',U')\).
\end{enumerate}
\end{prop}

\begin{proof}
Suppose \(\psi\) is a  \wcdns operator on \(X\) and let \(L\) be
a regular \gd set contained in an open set \(U\).  Define
\[\varphi(L,U) =
U\cap\bigcup \{\psi(J,W) \colon J \sq L, \mbox{ \(J\) is regular }
G_{\delta}, \mbox{ \(W\) is open}, \mbox{ } J \sq W \sq U\}.\]
Then \(\varphi(L,U)\) is open and \(L \sq \varphi(L,U) \sq U\)
and clearly \(\phi(L,U) \sq\phi(L',U')\) whenever
\(L \sq L'\) and \(U \sq U'\)

It remains to verify that $\phi$ is, indeed, a \wcdns operator.
So suppose that \(\varphi(L,U) \cap \varphi(K,V) \neq \n\). Then for some
regular \gd sets $L'$ and $K'$, and open sets $U'$ and $V'$, such that
$L'\sq L$, $K'\sq K$, $L'\sq U'\sq U$ and $K'\sq V'\sq V$, we have
$\psi(L',U')\cap\psi(K',V')\neq\n$. Hence either $\n\neq
L'\cap V'\sq L\cap V$ or $\n\neq K'\cap U'\sq K\cap U$, as required.
\end{proof}

On the other hand, it is not clear whether \cdns implies \mcdn.

In light of Theorem \ref{mngr2}, we might expect there to be a relationship
between m\(\delta\)n, \wcdns and \cdn.
Indeed, we have the following theorem.

\begin{thm} \label{ul1}
Each of the following properties of a space \(X\) implies the next:
\begin{enumerate}
\item \label{onei} Monotonically normal,
\item \label{twoii} m\(\delta\)n,
\item \label{threeiii} \mcdn,
\item \label{fouriv} \cdn,
\item \label{fivev} \wcdn,
\item \label{sixvi} m\(\delta\delta\)n.
\end{enumerate}
Moreover, every \mcdns space is
lm\(\delta\)n and every lm\(\delta\)n space is
m\(\delta\delta\)n.
\end{thm}

\begin{proof}
The proofs of (\ref{onei}) \(\rightarrow\) (\ref{twoii}),
(\ref{threeiii}) \(\rightarrow\) (\ref{fouriv}), (\ref{fouriv})
\(\rightarrow\) (\ref{fivev}) and the fact that
lm\(\delta\)n implies m$\dlt\dlt$n are trivial.

(\ref{twoii}) \(\rightarrow\) (\ref{threeiii}): We modify the proof
of Theorem \ref{mngr2}. Suppose \(H\) is an m\(\delta\)n operator
for \(X\) with \(H(L,K) \cap H(K,L) = \n\). Let \(L\) be a regular
\gd set and \(U\) an open set such that \(L \sq U\) and define
\(\psi(L,U) = H(L, X \sm U)\).  Then \(L \sq \psi(L,U) \sq
\ol{\psi(L,U)} \sq U\).  Assume \(L \cap V = \n\) and \(K \cap U =
\n\) where \(K\) is a regular \gd set contained in an open set
\(V\).  Then \(L \sq X \sm V\) and \(K \sq X \sm U\).  So by
monotonicity, \(\psi(L,U) \sq H(L,K)\).  Similarly, \(\psi(K,V) \sq
H(K,L)\). Therefore \(\psi(L,U) \cap \psi(K,V) = \n\). Monotonicity
of the operator \(\psi\) follows from the monotonicity of \(H\),
hence
\(\psi\) is a \mcdns operator for \(X\).

(\ref{fivev}) \(\rightarrow\) (\ref{sixvi}):
Again we modify the proof of Theorem \ref{mngr2}.
Suppose \(\psi\) is a  \wcdns operator for \(X\) and let \(L\)
and \(K\) be disjoint regular \gd sets in \(X\).  Define
\[H(L,K) = \bigcup \{\psi(J,U) \colon J \sq L \cap U, \mbox{ } J \mbox{ is
regular \(G_{\delta}\)}, \mbox{ } U \mbox{ is open, } U \cap K =
\n\}.\]
Then \(H(L,K)\) is open with \(L \sq H(L,K)\).  We  show that \(\ol{H(L,K)} \sq
X \sm K\).
Since \(X\) is \wcdn, if \(U\) is open with \(U \cap K = \n\) and
\(J\) is any regular \gd set contained in \(L \cap U\), then \(\psi(K, X \sm L)
\cap \psi(J,U) = \n\).
Hence \(\psi(K, X \sm L) \cap H(L,K) = \n\)
and so \(K \cap \ol{H(L,K)} = \n\).
 It is routine to show that the operator \(H\) is monotone.

To see that \mcdns implies lm\(\delta\)n,
assume \(\psi\) is a \mcdns operator for
\(X\). Let \(C\) and \(D\) be disjoint closed sets, \(C\) a regular
\gd set. Define \(H(C,D) = \psi(C,X \sm D)\). Then \(C \sq H(C,D)
\sq \ol{H(C,D)} \sq X \sm D\).  Suppose \(C \sq C'\) and \(D' \sq
D\).  Then \(X \sm D \sq X \sm D'\), hence \(H(C,D) \sq H(C',D')\).
\end{proof}

The proof of the following is routine.

\begin{prop}
Let \(M\) be a compact metrizable space. If $X\times M$ satisfies
any of the properties listed in Theorem \ref{ul1}, then so does $X$.
\end{prop}

\section{Factorizations of monotone normality}

Kohli and Singh \cite{sigma} factorize normality in terms of various
weak normality properties.  They define a space to be \(\Sigma\)-normal
if for each closed set \(C\) contained in an open set \(U\), there
exists a set \(W\) that is the complement of a regular \gd set such
that \(C \sq W \sq U\) and show that a space is
normal iff it is both weakly \(\delta\)-normal and \(\Sigma\)-normal.
There is an obvious monotone version of this result that factorizes
monotone normality into monotone $\Sigma$-normality and
m\(\delta\delta\)n. However, it turns out that we can do better than this in
the monotone case.

\begin{defn}
A space \(X\) is monotonically \(\Sigma\)-normal, or m$\Sigma$n, iff there is an
operator \(W\) assigning to each closed set $C$ and each open set
\(U\) containing $C$, an open set \(W(C,U)\) such that
\begin{enumerate}
\item $X \sm W(C,U)$ is a regular $G_\dlt$-set,
\item \(C \sq W(C,U) \sq U\) and
\item if \(C \sq C'\) and \(U \sq U'\), then \(W(C,U) \sq W(C',U')\).
\end{enumerate}
\end{defn}

\begin{prop} \label{prop1}
\(X\) is m$\Sigma$n iff there are operators
\(D\) and \(W\) assigning to each closed set $C$ and open set $U$,
containing $C$, sets $D(C,U)$ and $W(C,U)$ such that
\begin{enumerate}
\item $D(C,U)$ and $X \sm W(C,U)$ are regular \gd sets
\item \(C \sq D(C,U) \sq W(C,U) \sq U\),
\item $D(C,U)\cap W(X\sm U,X\sm C)=\n$,
\item if \(C \sq C'\) and \(U \sq U'\), and
then \(D(C,U) \sq D(C',U')\) and \(W(C,U) \sq W(C',U')\).
\end{enumerate}
\end{prop}

\begin{proof}
Suppose the conditions of the theorem hold, then clearly \(X\) is
m$\Sigma$n.  Conversely, suppose \(V\) is a
m$\Sigma$n operator for \(X\) and that \(C\sq
U\). Define $D'(C,U)=X \sm V(X\sm U,X\sm C)$, so that $C\sq D'(C,U)\sq U$ and
$D'(C,U)$ is a regular $G_\dlt$, and define $W(C,U)=V(D'(C,U),U)$.
It is routine to check conditions (1), (2) and (4). Now define
$D(C,U)=D'(C,U) \sm W(X\sm U,X\sm C)$, which is the intersection of two regular
$G_\dlt$-sets. Since
$W(X\sm U,X\sm C)=V\big(X\sm V(C, U),X\sm C\big)$ and $C\cap W(X\sm U,X\sm C)=\n$, we have operators $D$ and $W$ satisfying all four conditions.
\end{proof}

\begin{prop} \label{prop3} Every monotonically normal space and every
perfectly normal space is m$\Sigma$n.
\end{prop}

\begin{proof}
To show that every monotonically normal space is m$\Sigma$n, we extend the proof
that every normal space is \(\Sigma\)-normal \cite{sigma} and use the monotone
version of Urysohn's lemma \cite{stares}.

Suppose \(X\) is perfectly normal.  Then every open set
is the complement of a regular \gd set and defining
\(W(C,U) = U\) shows that $X$ is m$\Sigma$n.
\end{proof}

It turns out that a weaker property (that might be termed monotone
$\Sigma$ Hausdorff) is all that is needed to factorize monotone
normality in terms of m\(\delta\delta\)n.

\begin{defn}
A space \(X\) has \textit{property \((\star)\)} iff there are operators
$D$ and $E$ assigning to every $x\in X$ and
open set \(U\) containing $x$, disjoint sets \(D(x,U)\) and
\(E(x,U)\) such that
\begin{enumerate}
\item \(D(x,U)\) and \(E(x,U)\) are regular $G_\dlt$-sets,
\item \(x \in D(x,U) \sq U\) and
\item for every open set \(V\) and \(y \in V\), if \(x \notin V\) and \(y
\notin U\), then \(D(y,V) \sq E(x,U)\).
\end{enumerate}
\end{defn}

Of course, if $X$ is a regular space we can, without loss of generality,
drop the assumption that $D(x,U)\sq U$.

\begin{prop}
A space \(X\) has \textit{property \((\star)\)} iff there are operators
$D$ and $W$ assigning to each $x\in X$ and each open $U$ containing $x$,
sets \(D(x,U)\) and \(W(x,U)\) such that
\begin{enumerate}
\item $D(x,U)$ and $X \sm W(x,U)$ are regular $G_\dlt$-sets,
\item \(x \in D(x,U) \sq W(x,U) \sq U\) and
\item for every open set \(V\) and \(y \in V\), if \(x \notin V\) and \(y
\notin U\), then \(D(y,V) \cap W(x,U)= \n\).
\end{enumerate}
\end{prop}

\begin{proof} If $D$ and $E$ witness that \(X\) has property \((\star)\), define
$W(x,U)=X \sm E(x,U)$ for each $x\in U$.
If \(z \notin U\) and \(\hat{V} = X \sm D(x,U)\), then \(x \notin \hat{V}\)
and so \(z \in D(z,\hat{V}) \sq E(x,U)\).  Hence \(X \sm U \sq E(x,U)\) and
so \(W(x,U) \sq U\).  Since \(D(x,U)\) and \(E(x,U)\) are
disjoint, \(D(x,U) \sq W(x,U)\) and condition (3) is clear. The converse
follows just as easily.
\end{proof}

Property $(\star)$ is relatively easy to achieve.

\begin{thm} \label{propstar} Every m$\Sigma$n space
and every Tychonoff space with $G_\dlt$ points has property
$(\star)$.

Hence every monotonically normal space, every perfectly normal
space, every first countable Tychonoff space and every Tychonoff
space with a $G_\dlt$-diagonal has property $(\star)$.
\end{thm}

\begin{proof} Suppose that $X$ is m$\Sigma$n. Let $D$ and $W$ satisfy the
conditions of Proposition \ref{prop1}. Suppose that $U$ and $V$are
open sets and that $x\in U\sm V$ and $y\in V\sm U$. By (4),
$D(\{y\},V)\cap W(\{x\},U)\sq D(\{y\},X\sm\{x\}\cap
W(\{x\},X \sm \{y\})$, which is empty by (3). Hence $D(\{x\},U)$ and
$W(\{x\},U)$ define operators satisfying property $(\star)$.

Suppose now that $X$ is Tychonoff and has $G_\dlt$ points. Let $x\in U$.
Since $\{x\}$ is a $G_\dlt$-set, regularity implies that it is
a regular \gd set. Since $X$ is Tychonoff,
there is a continuous function \(f \colon X \to [0,1]\) such that \(f(x)=1\)
and \(f(X \sm U) = 0\). Define \(D(x,U) = \{x\}\) and
\(E(x,U) =f^{-1}(0)\). Then $D(x,U) $ and $E(x,U)$ are disjoint
regular $G_\dlt$-sets such that $x\in D(x,U)\sq U$ and $X \sm U\sq E(x,U)$, so that $D(y,V)\sq E(x,U)$, whenever $y\in V\sm U$.
\end{proof}

\begin{exmp} Assuming $\clubsuit^*$,
there is a space with property $(\star)$ that is not m$\Sigma$n.
\end{exmp}

\begin{proof} $\clubsuit^*$ asserts the existence of a sequence
$R_\alp=\{\bt_{\alp,n} \colon n\in\w\}$
for every limit ordinal $\alp\in\w_1$ that is cofinal in $\alp$ such
that, whenever $X$ is an uncountable subset of $\w_1$,
$\{\alp\in\w_1 \colon X\cap R_\alp\text{ is cofinal in }\alp\}$ contains a
closed unbounded set. $\clubsuit^*$ holds, for example, in any model of $V=L$.

Let $X=\w_1\times2$. For each limit $\alp$ and $n\in\w$,
define $B(\alp,n)=\{(\alp,1)\}\cup\{(\bt_{\alp,k},0) \colon n\les k\}$. Let $\T$
be the topology on $X$ generated by the collection
$$
\big\{\{(\alp,i)\} \colon \alp\text{ is a successor or }i=0\big\}
\cup\big\{B(\alp,n) \colon \alp\text{ is a limit}, n\in\w\big\}.
$$
With this topology, $X$ is zero-dimensional, hence Tychonoff, and
first countable, so has property $(\star)$.

If $U$ is an open set containing an uncountable subset of
$\w_1\times\{0\}$, for closed unboundedly many $\alp$,
$R_\alp\cap\{\bt \colon (\bt,0)\in U\}$
is cofinal in $\alp$, so that $\{\alp \colon (\alp,1)\in\ol{U}\}$ contains
a closed unbounded subset. Since the intersection of countably many closed
unbounded
subsets of $\w_1$ is, again, closed and unbounded, it follows that every
uncountable
regular \gd set in $X$ contains a closed unbounded subset of $\w_1\times\{1\}$.
Hence, if $C$ is any uncountable, co-uncountable subset of
$\w_1\times\{1\}$, $U=C\cup \left(\w_1\times\{0\}\right)$ and $D$ is any regular $G_\dlt$-set
containing $C$, then $C\sq U$, $U$ is open but $D\not\sq U$. Hence $X$ is not
m$\Sigma$n.
\end{proof}

Interestingly, property \((\star)\) is enough to push m\(\delta\delta\)n up to
monotone
normality.

\begin{thm} \label{starequivs} A space is monotonically normal iff it
has property $(\star)$ and is m\(\delta\delta\)n.
\end{thm}

\begin{proof}
Suppose \(H\) is an m\(\delta\delta\)n operator for \(X\) such that \(H(E,F) \cap
H(F,E) = \n\). Let \(U\) be an open set with \(x \in U\).  By property
\((\star)\),
there exist disjoint regular \gd sets \(D(x,U)\) and \(E(x,U)\) such that \(x
\in D(x,U) \sq U\) and for any open set \(V\) with \(x \notin V\), if \(y \in V
\sm U\) then \(D(y,V) \sq E(x,U)\).

Define \(\psi(x,U) = H(D(x,U), E(x,U))\).  Then \(D(x,U) \sq \psi(x,U)\), so \(x
\in \psi(x,U)\).  Suppose \(x \notin V\) and \(y \in V \sm U\).  Then by
monotonicity of \(H\), \(H(D(y,V), E(y,V)) \sq H(E(x,U),D(x,U))\).  It follows
that \(H(D(y,V),E(y,V)) \cap H(D(x,U), E(x,U)) = \n\).  Hence \(\psi(y,V)
\cap \psi(x,U) = \n\).  By Theorem \ref{mngr2}, \(X\) is monotonically
normal.

The converse is trivial given Theorems \ref{ul1} and \ref{propstar}.
\end{proof}

Hence, in any space with property $(\star)$, for example in a first
countable Tychonoff space, each of the properties listed in Theorem
\ref{ul1} is equivalent to monotone normality.

\begin{thm} \label{corr20}
\begin{enumerate}
\item If every point of $X$ is a regular \gd set, then
$X$ is monotonically normal iff it is \wcdn.
\item $X$ is \cdns iff it is  \wcdns and
\(\delta\)-normal.
\item If $X$ is normal, then $X$ is \cdns iff it is m\(\delta\delta\)n.
\end{enumerate}
\end{thm}

\begin{proof} In each case one implication follows from Theorem \ref{ul1}
and from the fact that a \cdns space is obviously
$\dlt$-normal.

To complete (1) and (2), suppose that $\psi$ satisfies conditions (1) and (2)
of Definition \ref{defnuln}. If every $x\in X$ is a regular $G_\dlt$, then
$\phi(x,U)=\psi(\{x\},U)$ satisfies conditions
(2) of Theorem \ref{mngr2} and $X$ is monotonically normal.
If $X$ is $\dlt$-normal and $L$ is a regular $G_\dlt$-subset of the open set
$U$, then there is an open set $\phi(L,U)$ such that
$L\subseteq\phi(L,U) \sq \ol{\phi(L,U)} \sq \psi(L,U) \sq U\). It is trivial to
check that, in this case, $\phi$ is a \cdns operator.

To complete (3),
suppose \(H\) is an m\(\delta\delta\)n operator for \(X\) with \(H(L,K) \cap
H(K,L) =
\n\).  Let \(L\) be a regular \gd set and \(U\) an open set such that
\(L \sq U\).  Since \(X\) is normal, there exists 
an open
set \(W_L\) such that \(W_L\) is the complement of a regular \gd set and \(L \sq
W_L \sq U\).  Define \(\psi(L,U) = H(L, X \sm W_L)\), then \(L \sq \psi(L,U) \sq
\ol{\psi(L,U)} \sq W_L \sq U\).  Now suppose \(L \cap V = \n\) and \(K
\cap U = \n\) where \(K\) is a regular \gd set contained in an open set
\(V\).  Then \(L \sq X \sm W_K\) and \(K \sq X \sm W_L\).
By monotonicity, \(\psi(L,U) \sq H(L,K)\) and \(\psi(K,V) \sq H(K,L)\), hence
\(\psi(L,U) \cap \psi(K,V) = \n\).
 Therefore \(\psi\) is a \cdns operator for \(X\).

\end{proof}

\section{Products with compact metrizable spaces and stratifiability}

A space $X$ is semi-stratifiable if there is an operator $U$
assigning to each $n\in \omega$ and closed set $D$ an open set $U(n,D)$
containing $D$ such that $\bigcap_{n \in \omega} U(n,D)=D$ and $U(n,D')\subseteq
U(n,D)$ whenever $D'\subseteq D$. If, in addition,
$\bigcap_{n \in \omega}\overline{U(n,D)}=D$, then $X$ is said to be stratifiable. A space $X$ is stratifiable iff $X\times M$ is monotonically normal for any (or all) infinite compact metrizable $M$ iff $X$ is both
semi-stratifiable and monotonically normal (see \cite{mn}).

\begin{defn}
A space \(X\) is \(\delta\)\textit{-semi-stratifiable} iff there is
an operator \(U\) assigning to each \(n \in \omega\) and regular \gd
set \(D\) in \(X\), an open set \(U(n,D)\) containing \(D\) such
that
\begin{enumerate}
\item[(1)] if \(E \subseteq D\), then \(U(n, E) \subseteq U(n, D)\) for each
\(n \in \omega\) and
\item[(2)] \(D = \bigcap _{n \in \omega} U(n,D)\).
\end{enumerate}
If in addition,
\begin{enumerate}
\item[(3)] \(D = \bigcap _{n \in \omega} \ol{U(n,D)}\),
\end{enumerate}
then \(X\) is \(\delta\)\textit{-stratifiable}.
\end{defn}

Just as for stratifiability, we may assume that the operator \(U\)
is also monotonic with respect to \(n\), so that \(U(n+1,D) \sq
U(n,D)\) for each \(n\) and regular \gd set \(D\).

The proof of the following is essentially the same as the proof of
the corresponding results for stratifiability and monotone
normality.

\begin{thm} \label{dsmddn} \label{dstrat3}
\begin{enumerate}
\item If \(X\) is
\(\delta\)-stratifiable, then $X$ is $\dlt$-semi-stratifiable and
m\(\delta\delta\)n.
\item If \(X\) is \(\delta\)-semi-stratifiable and lm\(\delta\)n, then it is
\(\delta\)-stratifiable.
\end{enumerate}
\end{thm}

\begin{thm}\label{dstratwo}
Let \(M\) be any infinite compact metrizable space. \(X\) is
\(\delta\)-stratifiable iff \(X \times M\) is
\(\delta\)-stratifiable iff $X\times M$ is $\mddn$.
\end{thm}

\begin{proof} Let \(\pi \colon X \times M \to X\) be the projection map. Since $M$
is compact, $\pi$ is both open and closed.

Suppose $X\times M$ is $\dlt$-stratifiable with
\(\delta\)-stratifiability operator $W$. By Theorem \ref{dsmddn},
$X\times M$ is $\mddn$.  To see that $X$ is $\dlt$-stratifiable, let
$D$ be a regular $G_\dlt$-subset of $X$. Fix some $r\in M$ and
define \(U(n,D) = \pi(W(n,D \times \{r\}))\). It is routine to
verify that $U$ is a $\dlt$-stratifiability operator for $X$.

Now suppose that $X$ is $\dlt$-stratifiable with operator $U$ such
that \(U(n, \n) = \n\) and satisfying \(U(n+1,E) \sq U(n,E)\) for
each \(n\) and regular \gd set \(E\). Suppose \(D\) is a regular
\(G_{\delta}\)-set in \(X \times M\).  Then \(D = \bigcap _{i \in
\omega} \ol{U}_i\) where \(D \sq U_i\) and \(U_i\) is open in \(X
\times M\) for each \(i\). Define \(D_r = D \cap (X \times \{r\})\)
for each \(r \in M\).  Then each \(D_r\) is a regular \gd set since
\(D_r = \bigcap _{i \in \omega} \ol{U_i \cap (X \times
B_{1/2^i}(r))}\) and \(D_r \sq U_i \cap (X \times B_{1/2^i}(r))\)
for all \(i \in \omega\). Clearly \(D = \bigcup _{r \in M} D_r\).
Moreover \(\pi(D_r)\) is a regular \gd set in \(X\) for each \(r \in
M\).

For each \(n \in \omega\) define
\[H(n,D) = \bigcup_{r \in M} U(n, \pi(D_r)) \times B_{\frac{1}{2^n}}(r).\]
We show that \(H\) is a \(\delta\)-stratifiability operator for \(X \times
M\).  Clearly \(H(n,D)\) is open for each regular \gd set \(D\) and \(n \in
\omega\).  That \(H\) is monotone is clear from the monotonicity of \(U\).
It is easily seen that \(D \sq H(n,D)\) for each \(n \in \omega\), so it remains
to prove that \(\bigcap_{n  \in \omega} \ol{H(n,D)} \sq D\).

Suppose \((x,s) \in \bigcap_{n \in \omega} \ol{H(n,D)}\sm D\).  Then there
exists
a basic open set \(V \ni x\) and \(k \in \omega\) such that \((V \times
B_{1/2^k}(s)) \cap D = \n\) and so \((V \times B_{1/2^k}(s)) \cap
(\pi(D_r) \times \{r\}) = \n\) for all \(r \in B_{1/2^k}(s)\).
Since \((x,s) \in \ol{H(n,D)}\) for each \(n \in \omega\), we may consider the
following two cases:

Case 1:  Assume \((x,s) \in \ol{\bigcup_{r \in B_{1/2^k}(s)} U(n, \pi(D_r))
\times B_{1/2^n}(r)}\) for all \(n \geqslant k+1\).  Then for all such \(n\),
\((W \times B_{1/2^m}(s)) \cap \bigcup_{r \in B_{1/2^k}(s)} U(n, \pi(D_r))
\times B_{1/2^n}(r) \neq \n\) for all basic open sets \(W \ni x\), \(m
\in \omega\).
It follows that for some \(t \in B_{1/2^k}(s)\), \(V \cap U(n, \pi(D_t)) \neq
\n\) for each \(n \geqslant k+1\).  Then, since \(U\) is monotonic with
respect to \(n\), \(V \cap \bigcap_{n \in \omega} U(n, \pi(D_t)) \neq
\n\).  Therefore \(V \cap \pi(D_t) \neq \n\), a contradiction.

Case 2:  Assume \((x,s) \in \ol{\bigcup_{r \notin B_{1/2^k}(s)} U(n, \pi(D_r))
\times B_{1/2^n}(r)}\) for all \(n \geqslant k+1\). Then for some \(p \notin
B_{1/2^k}(s)\), \((W \times B_{1/2^m}(s)) \cap (U(n, \pi(D_p)) \times
B_{1/2^n}(p)) \neq \n\) for all basic open sets \(W \ni x\), \(m \in
\omega\) and \(n \geqslant k+1\).  Thus, for all such \(m\) and \(n\),
\(B_{1/2^m}(s) \cap B_{1/2^n}(p) \neq \n\).  However, \(B_{1/2^{k+1}}(s)
\cap B_{1/2^n}(p) = \n\) for all \(n \geqslant k+1\), a contradiction.

Therefore \(D = \bigcap_{n \in \omega} \ol{H(n,D)}\) as required.

To complete the proof we wish to show that if $X\times M$ $\mddn$,
then $X$ is $\dlt$-stratifiable. Note first that we may assume that
$X\times\Om$ is $\mddn$, where $\Om=\w+1$ is the convergent
sequence. To see this note that if $W$ is a subspace of $M$ that is
homeomorphic to $\Om$, then any regular $G_\dlt$-subset of $X\times
W$ is in fact a regular $G_\dlt$-subset of $X\times M$, so that
$X\times W$ is also $\mddn$. The proof is now familiar.

Let $H$ be an $\mddn$ operator for $X\times\Om$ such that
$H(C,D)\cap H(D,C)=\nowt$ for any regular $G_\dlt$-sets $C$ and $D$.
For each $n\in\w$, let $\Om_n=(\w+1)\sm\{n\}$ and let $\pi \colon X\times
\Om\to X$ be the projection map. If $E$ is a regular $G_\dlt$-subset
of $X$ define
$$U(n,E)=\pi\big(H(E\times \{n\},X\times\Om_n)\big).$$

Clearly \(E \sq U(n,E)\) for each \(n\).
Suppose that $z\in \bigcap_{n\in \omega}\ol{U(n,E)}\sm E$. Then, as $E$ is
closed, there is a regular $G_\dlt$-set $D$ such that $z\in
D\subseteq X\sm E$. Hence $K=D\cap \bigcap_{n\in \omega}\ol{U(n,E)}$ is a
regular $G_\dlt$ such that $z\in K$, $K\cap E=\nowt$ and
$K\subseteq\bigcap_{n\in \omega}\ol{U(n,E)}$, from which it follows that
$$K\times\{w\}\subseteq \ol{\bigcup_{n\in\w}\ol{H(E\times\{n\},X\times\Omega_n)}}
=\ol{\bigcup_{n\in\w} H(E\times\{n\},X\times\Om_n)}.$$ 
Therefore, for some
$n\in \w$, we have
$$\nowt\neq H\big(K\times\{\w\},E\times\Om\big)\cap
H\big(E\times\{n\},X\times\Om_n\big),$$ but, by monotonicity, this
implies that
$$\nowt\neq H\big(K\times\{\w\},E\times\Om\big)\cap
H\big(E\times\Om,K\times\{\w\}\big),$$ which is a contradiction and
it follows that $\bigcap_{n\in \omega}\ol{U(n,E)}=E$.

\end{proof}

Clearly property \((\star)\) will have an effect on $\dlt$-stratifiability
although it not clear that it is productive.
Obviously, by Theorem \ref{propstar}, if
\(X\) and \(Y\) are Tychonoff with $G_\dlt$ points, in particular
if $Y$ is a compact metrizable space,
then \(X \times Y\) has property \((\star)\).
Furthermore, if the product of a space with some
compact metrizable space does not have property \((\star)\), then the space is
not stratifiable.

\begin{cor} \label{stdststar}
Let \(M\) be any infinite compact metrizable space. If $X\times M$
has property $(\star)$, in particular if $X$ is a Tychonoff space
with $G_\dlt$ points, then $X$ is stratifiable iff $X$ is
$\dlt$-stratifiable iff $X\times M$ is m\(\delta\delta\)n.
\end{cor}

\section{Examples}

The following lemma gives some simple sufficient conditions on the regular \gd
subsets of a space for it to be \wcdns or \mcdn.

\begin{lem} Let $X$ be a space.
\begin{enumerate}
\item If, whenever \(L\) and \(K\) are disjoint
regular \gd subsets, at least one of them is clopen, then $X$
is  \wcdn.
\item If every regular \gd subset of $X$ is clopen, then $X$ is both
\mcdns and $\dlt$-stratifiable.
\end{enumerate}\label{makes it ul}
\end{lem}

\begin{proof}
(1) For any regular \gd set \(L\) contained in an open set \(U\),
define \(\psi\) as follows:
\[\psi(L,U)=
  \begin{cases} L & \mbox{if \(L\) is clopen} \\
  U &  \mbox{if \(L\) is not clopen.}
  \end{cases}\]
Suppose \(L\) is clopen.  Then \(\psi(L,U) = L\) and \(\psi(K,V) \sq V\), where \(K\) is a regular \gd set contained in an open set \(V\).  Hence if
\(L \cap V = \n\) and \(K \cap U = \n\), then \(\psi(L,U) \cap \psi(K,V) = \n\).

(2) follows immediately by defining $\phi(L,U)=L$ and $U(n,L)=L$
for any $n\in\omega$ and regular \gd set $L$.
\end{proof}

Given a cardinal $\kp$, let $\mathbb L_\kp$ denote the space $\kp+1$
with the topology generated by isolating each $\alp\in \kp$ and
declaring basic open neighbourhoods of $\kp$ to take the form
$\mathbb L_\kp \sm C$, where $C$ is some countable subset of $\kp$. Note
that, if $\kp$ is uncountable, then any regular $G_\dlt$-subset of
$\mathbb L_\kp$ containing the point $\kp$ is clopen and
co-countable and that a regular $G_\dlt$-set that does not contain
$\kp$ is countable.

\begin{exmp} \label{Lindel1}
$\mathbb L_{\w_1}$ is monotonically normal and \(\delta\)-stratifiable, but
not semi-stratifiable. Moreover $\mathbb L_{\w_1}\times (\w+1)$ is
m\(\delta\delta\)n.
\end{exmp}

\begin{proof} By Lemma \ref{makes it ul} (2), $\mathbb L_{\w_1}$ is
$\dlt$-stratifiable. By Theorem \ref{mngr2}, defining \(\psi(x,U) =
U\), if $x=\w_1$, and \(\psi(x,U) = \{x\}\), otherwise, whenever $x$
is in the open set $U$, we see that $\mathbb L_{\w_1}$ is monotonically normal.
However, since $\{\w_1\}$ is not a \gd subset, $\mathbb L_{\w_1}$ is
not semi-stratifiable. That $\mathbb L_{\w_1}\times (\w+1)$ is
m\(\delta\delta\)n follows by Theorem \ref{dstratwo}.
\end{proof}

\begin{exmp} Let $\mathbb S$ be the Sorgenfrey line.
$\mathbb S$ is monotonically normal but not $\dlt$-stratifiable and
$\mathbb S\times(\w+1)$ is not m\(\delta\delta\)n.
\end{exmp}

\begin{proof} Since $\mathbb S\times(\w+1)$ is first countable and Tychonoff, it has
property $(\star)$. Since $\mathbb S$ is not stratifiable, $\mathbb
S\times(\w+1)$ is not monotonically normal and therefore not
m\(\delta\delta\)n.
\end{proof}

\begin{exmp} \label{Lindomega1}
$X=\big[\mathbb L_{\w_1}\times(\w+1)\big] \sm \{(\w_1,\w)\}$ is
\wcdn, but neither \cdns nor l$\mdn$.
\end{exmp}

\begin{proof} Let $T=\{(\alp,\w) \colon \alp\in\w_1\}$
and $R=\{(\w_1,k) \colon k\in\w\}$

To see that $X$ is not \cdn, note that $T$ is a regular $G_\dlt$-set
and that $U=X \sm R$ is an open set containing $T$. If $\phi(T,U)$
is any open set such that $T\sq \phi(T,U)\sq X \sm R$, then, for
some $k\in \w$, $\{(\alp,k) \colon (\alp,k)\in\phi(T,U)\}$ is
uncountable, so that $(\w_1,k)\in\ol{\phi(T,U)}$, but
$(\w_1,k)\notin U$. The same argument shows that $X$ is not l$\mdn$
either.

To see that $X$ is \wcdn, let \(L\) be a regular \gd subset
of the open set \(U\). First note that if $(\w_1,k)\in L$, then
$L\cap (\mathbb L_{\w_1}\times\{k\})$ is a clopen subset of $X$. For
each $(x,\w)\in L$, there is a least $k_x\in \w$ such that
$\{(x,j) \colon k_x\les j\}$ is a subset of $U$. Let
$B(x,U)=\{(x,\w)\}\cup\{(x,j) \colon k_x\les j\}$. Define
$$
\psi(L,U) = L \cup \bigcup \{B(x,U) \colon (x,\w)\in L\}.
$$
Then $L\subseteq \psi(L,U)\sq U$ and $\psi(L,U)$ is open.

Suppose that $L$ and $K$ are regular \gd sets, $U$ and $V$ are open
sets and that $L\sq U\sm V$ and $K\sq V\sm U$. Then
\begin{align*}\psi(L,U)& \cap \psi(K,V)\\ &=
\big(L \cup \bigcup \{B(x,U) \colon (x,\w)\in L\}\big)
\cap\big(K \cup \bigcup \{B(x,V) \colon (x,\w)\in K\}\big)\\
&=\bigcup \{B(x,U) \colon (x,\w)\in L\} \cap \bigcup \{B(x,V) \colon (x,\w)\in
K\}=\n,
\end{align*}
since otherwise, if $(x,k)\in\psi(L,U) \cap \psi(K,V)$, then
$(x,\w)\in L\cap K$.
\end{proof}

\begin{exmp} $X=\left[\lw\times\lww\right] \sm \{(\w_1,\w_2)\}$ is
\mcdns and $\dlt$-stratifiable, but not
$\mdn$.
\end{exmp}

\begin{proof} Let $L$ be a regular \gd subset of $X$ containing $(\w_1,\alp)$ (or
$(\alp,\w_2)$). Then $L$ contains  a clopen neighbourhood of
$(\w_1,\alp)$ (or $(\alp,\w_2)$). Hence every regular \gd subset of
$X$ is clopen and by Lemma \ref{makes it ul}, $X$ is \mcdns and $\dlt$-stratifiable.

To see that $X$ is not $\mdn$, suppose to the contrary that $H$ is an
$\mdn$ operator such that $H(C,D)\cap H(D,C)=\nowt$. For each
$\alp\in \w_1$ and $\bt\in \w_2$, let
\begin{align*}
&C_\alp=\big\{(\alp,\w_2)\big\},
&D_\alp=X \sm \big(\{\alp\}\times\lww\big),\\
&E_\bt=\big\{(\w_1,\bt)\big\}, &F_\bt=X \sm \big(\lw\times\{\bt\}\big).
\end{align*}
Notice that $C_\alp\cap D_\alp=E_\bt\cap F_\bt=\nowt$,
$C_\alp\subseteq F_\bt$, $E_\bt\subseteq D_\alp$,
$H(C_\alp,D_\alp)\subseteq\{\alp\}\times\lww$, and
$H(E_\bt,F_\bt)\subseteq\lw\times\{\bt\}$. Hence
$H(C_\alp,D_\alp)\subseteq H(F_\bt,E_\bt)$, so that
$H(C_\alp,D_\alp)\cap H(E_\bt,F_\bt)=\nowt$.

Now, for each $\bt \in \w_2$, there are no more than countably
$\alp\in\w_1$ such that $(\alp,\bt)\notin H(E_\bt,F_\bt)$. This
implies that there is a subset $W$ of $\w_2$ with cardinality $\w_2$
and some $\alp_0\in\w_1$ such that $(\alp_0,\w_1]\times\{\bt\}$ is a
subset of $H(E_\bt,F_\bt)$ for each $\bt\in W$. It follows that for
any $\alp_0\les \alp\in \w_1$ and any $\bt\in W$, $(\alp,\bt)\notin
H(C_\alp,D_\alp)$, so that $H(C_\alp,D_\alp)$ is not open, which is
the required contradiction.
\end{proof}

\end{document}